\def \version {May 2, 2021}

\documentclass[12pt]{article}
\usepackage{amsfonts,epsf,amsmath}
\usepackage{amssymb}
\usepackage{graphicx}
\usepackage{tikz}
\usepackage[cp1250]{inputenc}

\newtheorem{theorem}{\bf Theorem}
\newtheorem{corollary}[theorem]{\bf Corollary}

\newtheorem{proposition}[theorem]{\bf Proposition}
\newtheorem{problem}[theorem]{\bf Problem}
\newtheorem{conjecture}[theorem]{\bf Conjecture}

\newcommand{\qed}{\hfill $\Box$ \bigskip}

\begin{document}

\title{The $k$-path vertex cover: general bounds and chordal graphs}
\author{
Csilla Bujt{\'a}s $^{a,}$\footnote{The first author acknowledges the financial support from the  Slovenian Research Agency under the project N1-0108.}
\and
Marko Jakovac $^{b,c,}$\footnote{The second author acknowledges the financial support from the  Slovenian Research Agency (research  core  funding  No.\  P1-0297  and  research  projects  J1-9109, J1-1693 and N1-0095).}
\and
Zsolt Tuza $^{d,e,}$\footnote{Research supported in  part by the National Research, Development and Innovation Office -- NKFIH under the grant SNN 129364.}
}

\date{\small \version}
\maketitle

\vspace{-5mm}

\begin{center}
$^a$  Faculty of Mathematics and Physics, University of Ljubljana\\
Jadranska u.\ 19, 1000 Ljubljana, Slovenia\\
\texttt{csilla.bujtas@fmf.uni-lj.si}\\
\medskip
$^b$ Faculty of Natural Sciences and Mathematics, University of Maribor \\
Koro\v{s}ka cesta 160, 2000 Maribor, Slovenia \\
\medskip
$^c$ Institute of Mathematics, Physics, and Mechanics\\
Jadranska 19, 1000 Ljubljana, Slovenia\\
\texttt{marko.jakovac@um.si}\\
\medskip
$^d$ Alfr\'ed R\'enyi Institute of Mathematics\\
Re\'altanoda u.\ 13--15, 1053 Budapest, Hungary\\
\medskip
$^e$ Faculty of Information Technology, University of Pannonia\\
Egyetem u.\ 10, 8200 Veszpr\'em, Hungary\\
\texttt{tuza@dcs.uni-pannon.hu}
\end{center}


\begin{abstract}
For an integer $k\ge 3$, a $k$-path vertex cover of a graph $G=(V,E)$ is a set $T\subseteq V$ that shares a vertex with every path subgraph of order $k$ in $G$.
The minimum cardinality of a $k$-path vertex cover is denoted by  $\psi_k(G)$. We give estimates --- mostly upper bounds --- on  $\psi_k(G)$ in terms of various parameters, including vertex degrees and the number of vertices and edges. The problem is also considered on chordal graphs and planar graphs.
\end{abstract}

\bigskip



\section{Introduction}

A \emph{$k$-path vertex cover} (or \emph{$P_k$-transversal}) in a graph $G$  is a set $T\subseteq V(G)$ that shares at least one vertex with every path
 of order $k$ in $G$; i.e.,  $G-T$ does not contain any $P_k$ subgraph. 
In this paper we study the minimum cardinality of a $k$-path vertex cover in $G$,  denoted by $\psi_k(G)$. Remark that the case of $k=2$ means just the vertex cover (transversal) that is the complement of an independent set. 
It immediately follows that $\psi_2(G)=|V(G)|-\alpha(G)$, where $\alpha(G)$  denotes the size of a maximum stable set and is called the \emph{independence number} of $G$. 
Since these are central notions in graph theory and the literature contains many related deep results, here we concentrate on the cases with $k \ge 3$.
The case of $k=3$ is in close relation to the so-called dissociation number, which was firstly studied by Yannakakis \cite{ya-1981}. 
A subset of vertices in a graph $G$ is called a dissociation set if it induces a subgraph with maximum degree at most~$1$. 
The size of a maximum dissociation set in $G$ is called the \emph{dissociation number} of $G$ and is denoted by $\mathrm{diss}(G)$. Thus, $\psi_3(G)=|V(G)|-\mathrm{diss}(G)$. 
Several other connections can be found; for instance in \cite{Brause-2017},
  $k$-path vertex cover was shown to be related to $k$-path partition.

\paragraph{Application in Information Technology.}
The motivation for the invariant $\psi_k(G)$, which was introduced in \cite{Br-2011}, arises from communications in wireless sensor networks, where the data integrity is ensured by using Novotn\'y's $k$-generalized Canvas scheme \cite{novotny-2010}. 
The topology of wireless sensor networks can be represented by a graph, in which vertices represent sensor devices and edges represent communication channels between them.
We suppose that there are protected and unprotected sensor devices in the model. 
The attacker is unable to copy secrets from a protected device; however, an unprotected device can be captured by the attacker who can gain control over it. 
During the deployment and initialization of a sensor network, it should be ensured that at least one protected vertex exists on each path of order $k$ in the communication graph. 
The placement of protected sensors in a network is usually expensive. 
Thus an important goal is to minimize the cost of the network by minimizing the number of protected vertices, which coincides with the calculation of a minimum $k$-path vertex cover in a graph that represents such a network.

\paragraph{Known results on algorithmic complexity.}
It was proved in the pilot  paper that the computation of $\psi_k(G)$ is in general an NP-hard problem, but was shown to be solvable in linear time over the class of trees.

	The $k$-path vertex cover problem has gained much attention in the area of computational graph theory.
	Approximation algorithms for $\psi_3(G)$ were given by Tu et al.\ in \cite{tu-2011,tu-2013} and an exact algorithm for computing $\psi_3(G)$ 
	in running time\footnote{suppressing polynomially bounded factors by the
	 $O^*$-notation} $O^*(1.5171^n)$ 
	for a graph of order $n$ was presented in \cite{Kardos-2011}, and was later improved to $O^*(1.366^n)$ in \cite{xiko-2017}. The parametrized version  of the same problem was also considered. The goal is to decide whether there is a $k$-path vertex cover in $G$ of size at most $\ell$, where $k$ and $\ell$ are fixed positive integers. 
	A fixed parameter tractable (FPT) algorithm for the $3$-path vertex cover whose time complexity is $O^*(1.713^\ell)$ was given by Tsur in \cite{ts-2019}. Tsur recently presented a FPT algorithm for the $4$-path vertex cover which runs in $O^*(2.619^\ell)$ \cite{ts-2021}.
	 
	 The weighted version of $k$-path vertex cover was introduced by Bre\v sar et al.\ in \cite{Br-2014}.
In this version vertices are given weights and the problem is to find a minimum-weight set such that the graph obtained by deleting this set of vertices has no path $P_k$ as a subgraph.
Some special classes of graphs were considered, for instance, complete graphs and cycles, and an algorithm that computes the weighted $k$-path vertex cover number of a tree with time complexity $O(k \cdot |V(G)|)$ was presented. Approximation algorithms for the weighted and connected versions of the $k$-path vertex cover were considered in \cite{Lizh-2016}, and later also in \cite{Ran-2019} where the emphasis was primarily on $k=3$. In both cases an extra assumption was present, namely that the subgraph induced by a $k$-path vertex cover was connected. The connected version of the $k$-path vertex cover was also studied in \cite{Chen-2018, Lizh-2020}.

\paragraph{Known results on graph products.}
The $k$-path vertex cover problem was heavily studied on several types of graphs products.
Already in \cite{Br-2013} the Cartesian product of two paths, i.e.\ grid graphs, was considered.
In \cite{Li-2018} those results were extended to the Cartesian product of three paths.
The results on grid graphs were later improved in \cite{jata-2013}, and the ideas of the proofs were used on the strong product of paths.
In the same paper an upper an a lower bound were given for the $k$-path vertex cover of the lexicographic product of arbitrary graphs.
The upper bound was shown to be tight for any choice of factor graphs when $k=3$. The $k$-path vertex cover of the rooted product of arbitrary graphs was considered in \cite{jakovac-2015}.

\paragraph{Some known lower and upper bounds.}
Throughout many articles on the $k$-path vertex cover the following problem appeared to gain much attention.

\begin{problem}   \label{p:ab}
 For a given\/ $k$, determine the set of pairs\/ $(a,b)$ of nonnegative reals
  such that\/ $\psi_k(G)\le an+bm$ holds for every graph with\/ $n$ vertices
  and\/ $m$ edges.
 Study the same also for interesting classes of graphs.
\end{problem}

In the next theorems we cite some tight related results from the paper of
Bre\v sar et al.\ \cite{Br-2011}, which will be applied later.

\begin{theorem}[\cite{Br-2011}]  \label{p:tree}
	Let\/ $T$ be a tree of order\/ $n$ and\/ $k$  an integer with\/ $k\ge 3$. Then, we have\/
	$\psi_k(T) \le n/k$. 
\end{theorem}

\begin{theorem}[\cite{Br-2011}]
	\label{thm:br-2011}
	Let\/ $G$ be a graph of order\/ $n$ and size\/ $m$. Then
	$$\psi_3(G) \le \frac{2n+m}{6} \qquad \mbox{and} \qquad  \psi_3(G) \le \frac{m}{2}.$$
	 Moreover, if\/ $G$ is subcubic, then\/ $\psi_3(G) \le n/2$ also holds.
\end{theorem}

Note that $\psi_3(G) \le m/2$ follows from the procedure in which we sequentially remove a vertex of degree at least $2$ from $G$ and put it into the $3$-path vertex cover. 
During this procedure, we delete at most $m/2$ vertices and what remains at the end is a graph with maximum degree at most 1, hence $P_3$-free. Further, if $\Delta(G) \ge 3$, then the strict inequality $\psi_3(G) < m/2$ holds as we may remove a vertex and at least three incident edges in the first step of the procedure.

The relation $\psi_3(G) \le n/2$ can also be obtained via an iterated algorithmic process.
Here, starting with any partition $(V_1, V_2)$ of the vertex set of a subcubic graph, a vertex which has more than one neighbor
 in its partition class is moved to the other class, and  hence increasing the size of the cut, until no such vertices remain.
Then each of $V_1$ and $V_2$ is a 3-path vertex cover.

In~\cite{Br-2013} the upper bound
$$\psi_3(G) \leq \frac{\ell}{\ell+2}n + \frac{1}{(\ell+1)(\ell+2)}m$$
is established, where $\ell = \left\lceil \frac{m}{n} \right\rceil -1$.
Also this bound is tight.

\paragraph{Our results.}
In this paper we continue this track of research and prove upper bounds on $\psi_k(G)$ in terms of various parameters, including vertex degrees and the number of vertices and edges.
We generalize the result of \cite{Br-2013} for $d$-regular graphs,
  \begin{equation}   \label{eq:reglow}
  \psi_k(G) \geq \frac{d-k+2}{2d-k+2}n,
 \end{equation}
to arbitrary graphs in terms of minimum and maximum degree.
We also consider the $k$-path vertex cover problem on some famous classes of graphs, such as chordal graphs and planar graphs.
Our paper is structured as follows. In Section~\ref{sec:2} we study Problem~\ref{p:ab}  
and prove general upper bounds on $\psi_3(G)$, $\psi_4(G)$, and $\psi_k(G)$ in terms of the order and size of $G$.
 
We also present a method that can be applied to generate, from feasible pairs
 $(a,b)$ of Problem~\ref{p:ab}, further feasible pairs $(a',b')$.
In Section~\ref{sec:3} we give lower and upper bounds on $\psi_k(G)$ with functions of vertex degrees.
Chordal graphs are considered in Section~\ref{sec:4}. Section \ref{sec:5} describes a way to improve the estimates further, and contains more  open problems. In particular, Subsection~\ref{sec:planar} discusses the problem on planar graphs. Some further problems are mentioned in the concluding section.

\subsection*{Notation}

We use standard notation as follows. For a simple undirected graph $G$ we usually denote the order $|V(G)|$ and the size $|E(G)|$ by $n$ and $m$. 
The notations $\delta(G)$, $\Delta(G)$, $\chi(G)$, and $\omega(G)$ stand for the \emph{minimum vertex degree}, the \emph{maximum vertex degree}, the \emph{chromatic number}, and the \emph{clique number} of $G$, respectively. 
We will simply write $\delta$, $\Delta$, $\chi$, and $\omega$ if $G$ is clear from the context.
The \emph{degree} of a vertex $v \in V(G)$ in graph $G$ is denoted by $d_G(v)$, and the set of its neighbors is denoted by $N_G(v)$ and called the \emph{open neighborhood} of $v$. 
Moreover, we define $N_G[v]=N_G(v) \cup \{v\}$ as the \emph{closed neighborhood} of $v$.
Also here, if $G$ is understood, we will simply write $d(v)$, $N(v)$ and $N[v]$.
The average vertex degree is $\overline{d}=\overline{d}(G)=2m/n$.
The subgraph induced by vertex set  $Y\subseteq V(G)$ in $G$ is referred to as $G[Y]$.
Notations $P_k$, $C_k$, and $K_k$ stand for the  path, cycle, and complete graph of order $k$, respectively. 


\section{Upper bounds related to Problem~\ref{p:ab}} \label{sec:2}

Here our goal is to make a modest step in the direction of solving Problem \ref{p:ab}. 
Suppose that $an+bm$ is a general upper bound on $\psi_k(G)$ for an integer $k \ge 3$.
If $b=0$, the best estimation is obtained with $a=1$, as shown by complete graphs which have $\psi_k(K_n) = n-k+1$; and for $a=0$  with $k=3$ 
the best pair is $b=1/2$. The corresponding upper bound is attained by graphs in which every component is a $P_3$. Also, $(a,b)=(1/3,1/6)$ provides a valid upper bound, by Theorem~\ref{thm:br-2011}.

\begin{theorem} \label{thm:nm4}
	If\/ $G$ is a graph of order\/ $n$ and size\/ $m$, then
	$$\psi_3(G) \le \frac{n+m}{4},$$
 	where equality holds if and only if all components of $G$ are isomorphic to $C_4$.
\end{theorem}
\begin{proof}
	We proceed by induction on $n$ and suppose that $G$ is given on  $n \ge 3$ vertices. If $\Delta(G) \le 3$, the inequalities $\psi_3(G) \le n/2$ and $\psi_3(G) \le m/2$
	 from Theorem~\ref{thm:br-2011} imply the statement. 
	 Otherwise, there exists a vertex $v \in V(G)$ with $d_G(v) \ge 4$.
By the induction hypothesis, $G'=G-v$ satisfies the inequality. Note that $G'$ has $n'=n-1$ vertices and  $m'\le m-4$ edges and further, if $T'$ is a $3$-path vertex cover in $G'$, then $T' \cup \{v\}$ is a $3$-path vertex cover in $G$. We may conclude that
	$$\psi_3(G) \le \psi_3(G')+1 \le \frac{n'+m'}{4} +1 \le \frac{n+m-5}{4} +1 < \frac{n+m}{4}.$$
	This completes the proof for the upper bound.
	
	Now, we consider the graphs attaining the upper bound in the theorem.  Let $G$ be any graph with $\psi_3(G) = \frac{n+m}{4}$.
	It suffices to prove that if $G$ is connected, then $G \cong C_4$.
	It follows from Theorem~\ref{thm:br-2011} that $n=m$, for otherwise
	$\psi_3(G) \le \min\{\frac{n}{2},\frac{m}{2}\} < \frac{n+m}{4}$. As we discussed after Theorem~\ref{thm:br-2011}, $\Delta(G) \ge 3$ would imply $\psi_3(G) < m/2$ and, therefore, $\psi_3(G) < \frac{n+m}{4}$ which is a contradiction. Hence, $G$ satisfies $\Delta(G) \le 2$ and $n=m$. This implies $G \cong C_n$. 
	Then $\psi_3(G)=\left\lceil \frac{n}{3} \right\rceil$, which is easily seen to be smaller than $n/2$ unless $n=4$.
	\qed
\end{proof}

For the upper bound stated in Theorem~\ref{n-m-thm} below, we also have a sharp example, namely $H= K_6 -{\cal M}$ which is obtained from the complete graph $K_6$ by removing a perfect matching ${\cal M}$. Then we have $n=6$, $m=12$, and $\psi_3(H)=4$.

\begin{theorem}\label{n-m-thm}
	Let\/ $G$ be a graph of order\/ $n$ and size\/ $m$. Then
	$$\psi_3(G) \leq \frac{4n+m}{9}.$$
\end{theorem}
\begin{proof}
	By Theorem~\ref{thm:br-2011}, every graph satisfies 	$\psi_3(G) \leq \frac{2n+m}{6}$. To show that $\psi_3(G) \leq \frac{4n+m}{9}$ is also true, we first note that the inequality clearly holds for graphs of order at most $2$. Then, we proceed by induction on the order $n$ of $G$ and consider two cases.
	\begin{itemize}
		\item If $m \le 2n$, the inequality chain
		$$\frac{4n+m}{9} = \frac{n}{3}+\frac{m}{6}+\left(\frac{n}{9}-\frac{m}{18}\right) \ge \frac{n}{3}+\frac{m}{6} \ge \psi_3(G)$$
		establishes the statement without using the induction hypothesis.
		\item  If $m >2n$, then $\Delta(G) \ge 5$ follows. Let $v $ be a vertex of degree at least $5$. Deleting $v$ from $G$, we obtain the graph $G'$ which is of order $n'=n-1$ and of size $m' \le m-5$. By the induction hypothesis, $G'$ satisfies the inequality, i.e.\ $\psi_3(G')\le (4n'+m')/9$.
Moreover, if $T'$ is a $3$-path vertex cover in $G'$, then $T' \cup \{v\}$ covers all $P_3$ subgraphs in $G$. This yields $\psi_3(G) \le \psi_3(G') +1$. Then, we have the desired result by the following relations:
		$$\psi_3(G) \le \psi_3(G') +1 \le \frac{4n'+m'}{9} +1 \le \frac{4(n-1)+(m-5)}{9} +1= \frac{4n+m}{9}.$$	
	\end{itemize}
	This completes the proof.
	\qed	
\end{proof}

We also present a general result that shows how one can obtain new feasible pairs $(a',b')$ if a feasible pair $(a,b)$ is already given.

\begin{theorem}
	\label{thm:gen-n-m}
	Let\/ $a$ and\/ $b$ be nonnegative reals and\/ $k$ an integer with\/ $k\ge 3$ such that\/ $\psi_k(G) \le an+bm$ holds for every graph\/ $G$. Moreover, let\/ $x$ be an arbitrary positive integer which satisfies 
	\begin{equation}
	\label{eq:1}
	\frac{1-a-b}{b} \le x \le \frac{2-2a}{b},
	\end{equation}
	and define\/ $a'$ and\/ $b'$ as follows:
	$$a'=a+bx- \frac{x(1+b-a)}{x+2}, \qquad  \qquad b'=\frac{2-2a-bx}{x+2}.
	$$
	Then,  the inequality 
	\begin{equation}
	\label{eq:2}
	\psi_k(G) \le a'n+b'm
	\end{equation}
	holds for every graph\/ $G$.
\end{theorem}
\begin{proof}
	Let $a$, $b$, $a'$, $b'$ and $x$ satisfy the conditions in the theorem.
	We introduce the notation
	$$y= \frac{x(a+b+bx-1)}{x+2}.$$
	By (\ref{eq:1}), we have $a+b+bx-1 \ge 0$ that implies $y \ge 0$.
	One can check that $a'=a+y$ and $b'=b-\frac{2y}{x}$. As the second inequality in (\ref{eq:1}) is equivalent to $bx \ge 2y$, we have  $b'\ge 0$ and $a' \ge a \ge 0$. Consequently, the right-hand side of (\ref{eq:2}) is always nonnegative and inequality (\ref{eq:2}) clearly holds for small graphs with $n \le k-1$. We proceed by induction on $n$. Consider a graph $G$ of order $n$ and size $m$.
	\begin{itemize}
		\item First, suppose that  $m \le \frac{xn}{2}$. This implies $n - \frac{2m}{x} \ge 0$. Using this relation and rewriting the formula $a'n+b'm$, we get
		$$a'n+b'm= (a+y)n +\left(b-\frac{2y}{x} \right)m= (an+bm)+y\left(n-\frac{2m}{x}\right) \ge an+bm \ge \psi_k(G)$$
		that verifies (\ref{eq:2}) for the first case.
		\item Next, suppose that the strict inequality $m >\frac{xn}{2}$ holds. It implies $\bar{d}(G) >x$  and therefore
		 $\Delta(G) \ge x+1$ holds. Let $v$ be a vertex with $d_G(v) \ge x+1$. Deleting $v$ from $G$, we obtain a graph $G'$ with $n'=n-1$ vertices and  $m' \le m-(x+1)$ edges. Since every edge that was deleted from $G$ is covered by $v$, we have  $\psi_k(G) \le \psi_k(G') +1$.
		By our hypothesis, $G'$ satisfies (\ref{eq:2}). The following computation proves the same property for $G$:
		\begin{align*}
		a'n+b'm & \ge (a'n'+b'm')+ a' +b'(x+1)\\  
		& \ge \psi_k(G') + (a+y) + \left(b-\frac{2y}{x}\right)(x+1)\\
		& = \psi_k(G') +  a+bx+b -\frac{x+2}{x}\, y\\
		& = \psi_k(G') +  a+bx+b -(a+b+bx-1) \\
		& = \psi_k(G') +  1 \ge \psi_k(G).
		\end{align*} 
	\end{itemize}
This finishes the proof of the theorem. \qed
	\end{proof}

Since the formulation of the above theorem is quite technical, we present some upper bounds which can be obtained starting with Theorem~\ref{n-m-thm} and applying Theorem~\ref{thm:gen-n-m} iteratively. We begin with the case of $k=3$ and with the upper bound $\frac{4}{9}n+\frac{1}{9}m$ and set $a_0=4/9$, $b_0=1/9$. Then for each $1 \le i \le 10$, we apply Theorem~\ref{thm:gen-n-m} with $x=i+4$, $a=a_{i-1}$, and $b=b_{i-1}$.

\begin{corollary}
\label{many-bounds}
	For every graph\/ $G$ of order\/ $n$ and size\/ $m$, the following inequalities hold:
	\begin{align*}
	\psi_3(G) & \le \frac{11}{21} n + \frac{5}{63 }m,\\
	\psi_3(G) & \le \frac{7}{12} n + \frac{5}{84}m,\\
	\psi_3(G) & \le \frac{17}{27} n + \frac{5}{108}m,\\
	\psi_3(G) & \le \frac{2}{3} n + \frac{1}{27}m,\\
	\psi_3(G) & \le \frac{23}{33} n + \frac{1}{33}m,\\
	\psi_3(G) & \le \frac{13}{18} n + \frac{5}{198}m,\\
	\psi_3(G) & \le \frac{29}{39} n + \frac{5}{234}m,\\
	\psi_3(G) & \le \frac{16}{21} n + \frac{5}{273}m,\\
	\psi_3(G) & \le \frac{7}{9} n + \frac{1}{63}m,\\
	\psi_3(G) & \le \frac{19}{24} n + \frac{1}{72}m.
	\end{align*}
\end{corollary}


\begin{proposition} \label{Delta-2}
	If\/ $G$ is a graph of maximum degree\/ $\Delta(G) \le 2$, then 
	$$\psi_k(G) \le \frac{2}{k+1}\, n \qquad \mbox{and} \qquad \psi_k(G) \le \frac{2}{k+1}\, m $$
	hold for every\/ $k \ge 3$.
  Further, the first inequality is tight if and only if each component of\/ $G$ is a\/ $(k+1)$-cycle; the second one is tight if and only if each component is either an isolated vertex or a\/ $(k+1)$-cycle.

\end{proposition}
\begin{proof}
	It follows from $\Delta(G) \le 2$ that every component of $G$ is either a path or a cycle. Since $\psi_k$ is an additive invariant and 
	$$\psi_k(P_n)= \left\lfloor \frac{n}{k}\right\rfloor = \left\lfloor \frac{m(P_n)+1}{k}\right\rfloor \le 
	\psi_k(C_n)= \left\lceil \frac{n}{k}\right\rceil = \left\lceil \frac{m(C_n)}{k}\right\rceil $$
	 holds for every $n \ge k$, moreover $\psi_k(P_n)= \psi_k(C_n)=0$ holds if $ n <k$, 
	  we easily derive  
	 $$\psi_k(G) \le \frac{2}{k+1}\, m  \le \frac{2}{k+1}\, n $$
	 and identify the sharp cases.
	 \qed
\end{proof}

\begin{theorem} \label{thm:psi-4}
Let\/ $G$ be a graph of order\/ $n$ and size\/ $m$. Then
$$\psi_4(G) \leq \frac{n+3m}{10}.$$
\end{theorem}
\begin{proof}
 We proceed by induction on $n$.
 If $\Delta(G) \le 2$, then $5\psi_4(G) \le 2n$ and $15\psi_4 \le 6m$ follow by  Proposition~\ref{Delta-2}. Therefore, we have $20\psi_4(G) \le 2n +6m$ and the statement follows.
 If $\Delta(G) \ge 3$, there is a vertex $v$ of degree at least $3$ and we may consider the graph $G'=G-v$. Since $G'$ has $n'=n-1$ vertices, $m' \le m-3$ edges and, by our hypothesis, it satisfies $\psi_4(G') \leq (n'+3m')/10$,
 $$\psi_4(G) \le \psi_4(G') +1 \le \frac{n'+3m'}{10} +1 \le \frac{(n-1)+3(m-3)}{10} +1= \frac{n+3m}{10}$$
is obtained that completes the proof.
\qed
\end{proof}

Theorem~\ref{thm:psi-4} holds with equality if each component of $G$ is a $C_5$.
 (Currently the 5-cycle is the only connected graph for which we know the equality.)
Moreover, using Theorem~\ref{thm:gen-n-m}, further upper bounds on $\psi_4$ can be generated.
However, we do not list such bounds here.
Instead, generalizing Theorem~\ref{thm:gen-n-m} we close this section with describing a general class of recursions showing how further feasible pairs $(a',b')$ can be obtained from an already known feasible pair $(a,b)$.

\begin{theorem}
Assume that\/ $\psi_k(G) \le an+bm$ is a universally valid inequality for all graphs\/ $G$.
Suppose further that the real numbers\/ $q,w>0$ and the integer\/ $x>0$ satisfy the following inequalities:
 $$
   w\le 2q , \qquad a+qx < 1 , \qquad
    w \le \frac{(a+b-1)+(q+b)x}{x+1} .
 $$
Then\/ $\psi_k(G) \le a'n+b'm$ also holds for every graph\/ $G$ with the values
 $$
   a' = a+qx , \qquad b' = b-w .
 $$
\end{theorem}

\begin{proof}
Let us note first that under the given conditions we have $b'>0$.
Indeed, from the assumed inequalities we obtain
 $$
   1 \le a + qx + bx + b - (x+1)w < 1 + (b-w)(1+x) .
 $$
If $G$ has average degree at most $x$, then $m\le xn/2$ holds, thus
\begin{align*}
   \psi_k(G) &\le an + bm = (a'-qx) n + (b'+w) m = a'n + b'm - (qxn - wm)\\
						 &\le a'n + b'm - (qxn - 2qm) \le a'n + b'm .
\end{align*}
Otherwise there exists a vertex $v$ whose degree is at least $x+1$.
Supplementing any $k$-path vertex cover of the graph $G'=G-v$ with $v$
we obtain a $k$-path vertex cover of $G$.
Here $G'$ has order $n'=n-1$ and size $m'\le m-x-1$.
Hence, assuming by induction that $(a',b')$ is a feasible pair for the
smaller graph $G'$, we obtain:
\begin{align*}
   \psi_k(G) &\le 1 + \psi_k(G') \le 1 + a'(n-1) + b'(m-x-1)\\
						 &= a'n + b'm - ( a+qx + (b-w)(x+1) -1 )\\
						 &\le a'n + b'm .
\end{align*}
\qed
\end{proof}

Based on this theorem, one can derive an infinite class of feasible pairs
 $(a',b')$ already from one single pair $(a,b)$.
Nevertheless, we think that the most interesting case occurs when the two
 upper-bound conditions on $w$ coincide.
This special value is the one exposed in Theorem \ref{thm:gen-n-m} above.
More explicitly, it is obtained by substituting $q=\frac{(a+b+bx-1)}{x+2}$
 and $w=2q$.

\section{Estimates in terms of vertex degrees }
\label{sec:3}

In this section we give general estimations on $\psi_k(G)$
 in terms of the order and of the minimum, maximum, and average vertex degrees.

\subsection{In terms of $\delta$ and $\Delta$}

In the proof of the first general upper bound we
 will refer to the following decomposition theorem of Lov\' asz.

\begin{theorem}[\cite{L-1966}]  \label{Lovasz-decomp}
If\/ $a$ and\/ $b$ are nonnegative integers and\/ $G$ is a graph of maximum degree at most\/ $a+b+1$, then the vertex set of\/ $G$ can be partitioned into two sets which induce subgraphs of maximum degree at most\/ $a$ and\/ $b$.
\end{theorem}

The following theorem was proved by Bre\v sar \textit{et al.} \cite{Br-2011} in 2011. 
\begin{theorem}[\cite{Br-2011}] \label{thm:Br-Delta}
Let\/ $G$ be a graph of order\/ $n$ and of maximum degree\/ $\Delta$. Then
$$\psi_3(G) \le \frac{\left\lceil\frac{\Delta-1}{2}\right\rceil}{\left\lceil\frac{\Delta+1}{2}\right\rceil}\, n.$$
\end{theorem}

By definition, $\psi_k(G) \le \psi_3(G) $ if $k \ge 3$.  Theorem~\ref{thm:Br-Delta} therefore implies $\psi_k(G) \le \frac{\Delta-1}{\Delta+1}\,n$ for every $k \ge 3$ if $\Delta$ is odd and $\psi_k(G) \le \frac{\Delta}{\Delta+2}\,n$ if $\Delta$ is even. For $k=3$ and any odd $\Delta$ the bound is tight, as shown by the complete graphs  $K_{\Delta+1}$. Similarly, omitting a perfect matching from $K_{\Delta+2}$ we obtain tight examples for all even $\Delta$ with $k=3$.

Our next result shows that the upper bound can be improved for all $k \geq 4$ and all $\Delta \geq 4$. Note that part $(i)$ of Theorem~\ref{delta-thm} gives the same bound as  Theorem~\ref{thm:Br-Delta} for $k=3$ and that the bound stated for $\Delta=2$ here is exactly the same as the one in Proposition~\ref{Delta-2}.


\begin{theorem}\label{delta-thm}
Let\/ $k$ and\/ $\Delta$ be integers with\/ $k\ge 3$ and\/ $\Delta=2$ or\/ $\Delta \ge 4$,
 and let\/ $G$ be a graph of order\/ $n$ and of maximum degree at most\/ $\Delta$.
Then the following hold.
\begin{itemize}
  \item[$(i)$] If\/ $\Delta\geq 2$ is even, then
   $$\psi_k(G)\le \frac{(k-1)(\Delta-2)+4}{(k-1)\Delta+4}\,n.$$
  \item[$(ii)$] If\/ $\Delta\geq 5$ is odd, then
   $$\psi_k(G)\le \frac{(k-1)(\Delta-3)+8}{(k-1)(\Delta-1)+8}\,n.$$
\end{itemize}
\end{theorem}

\begin{proof}
We fix an integer $k\ge 3$ and proceed by induction on $\Delta$.
If $\Delta(G)\le 2$ then, by Proposition~\ref{Delta-2}, $\psi_k(G) \le \frac{2}{k+1}\, n$ holds.
This verifies the basic case of $(i)$.
Now, assume that $\Delta$ is even, $\Delta\geq 4$, and that $G$ has maximum degree at most $\Delta$.
By Theorem~\ref{Lovasz-decomp}, $V(G)$ can be partitioned into two sets $A$ and $B$,
 such that $\Delta(G[A]) \le \Delta-2$ and $\Delta(G[B]) \le 1$.
Since  $G[B]$ is $P_k$-free,  $\psi_k(G) \le |A|$ holds. On the other hand, if $T$ is a minimum $k$-path vertex cover in $G[A]$, then $B\cup T$ covers all paths of 
order $k$ in $G$, and we have $\psi_k(G) \le |B|+|T|$.
If $|A|\le \frac{(k-1)(\Delta-2)+4}{(k-1)\Delta+4}\,n$
holds, then $\psi_k(G) \le |A|$ establishes the upper bound. In the
other case, we  observe that
\begin{eqnarray}
 \psi_k(G) & \le & |B|+|T|
  \ \le \ n-|A|+\frac{(k-1)(\Delta-4)+4}{(k-1)(\Delta-2)+4}\,|A| \nonumber \\
    & = & n- \frac{2(k-1)}{(k-1)(\Delta-2)+4}\,|A|
     \ < \ \frac{(k-1)(\Delta-2)+4}{(k-1)\Delta+4}\,n. \nonumber
\end{eqnarray}
This proves $(i)$.

Suppose next that $\Delta$ is odd and $\Delta \ge 5$.
Say, $\Delta=2a+2b+1$ where $a,b$ are positive integers.
From Theorem~\ref{Lovasz-decomp} we know that there exists a vertex partition
 $V(G)=A\cup B$ such that $\Delta(G[A]) \le 2a$ and $\Delta(G[B]) \le 2b$.
Moreover, by $(i)$, there exist $k$-path vertex covers $T_A$ in $G[A]$ and
 $T_B$ in $G[B]$ with
  $$
    |T_A| \leq \frac{(k-1)(a-1)+2}{(k-1)a+2}\,|A| , \qquad
      |T_B| \leq \frac{(k-1)(b-1)+2}{(k-1)b+2}\,|B| .
  $$
Both $T_A\cup B$ and $T_B\cup A$ are $k$-path transversals of $G$.
Hence, denoting $x=|A|/n$ and $y=|B|/n=1-x$, we have
 $$
   \frac{\psi_k(G)}{n} \leq \min \left( \frac{(k-1)(a-1)+2}{(k-1)a+2}\,x + (1-x) ,
     \frac{(k-1)(b-1)+2}{(k-1)b+2}\,(1-x) + x \right).
 $$
The first term of minimization is a decreasing function of $x$ while the
 second term is an increasing function.
 Hence, the upper bound never exceeds the one obtained when the two numbers are equal, that is
\begin{align*}
   \frac{(k-1)(a-1)+2}{(k-1)a+2}\,x + (1-x) &= \frac{(k-1)(b-1)+2}{(k-1)b+2}\,(1-x) + x ,\tag{3}\\ 
   \frac{k-1}{(k-1)b+2} &= \left( \frac{k-1}{(k-1)a+2} + \frac{k-1}{(k-1)b+2} \right) x ,\\
   1 &= \left( \frac{(k-1)b+2}{(k-1)a+2} + 1 \right) x .
 \end{align*}
Recalling that $2a+2b=\Delta-1$, the worst case occurs when
 $$
   x = \frac{(k-1)(2a)+4}{(k-1)(\Delta-1)+8} , \qquad
      y = \frac{(k-1)(2b)+4}{(k-1)(\Delta-1)+8} .
 $$
Thus, substituting this particular value of $x$ into the left-hand side of $(3)$ we obtain
 \begin{align*}
   \frac{\psi_k(G)}{n} &\leq \frac{[(k-1)(2a-2)+4] + [(k-1)(2b)+4]}{(k-1)(\Delta-1)+8}\\
											 &=\frac{(k-1)(\Delta-3)+8}{(k-1)(\Delta-1)+8}. 
 \end{align*}
This proves $(ii)$ and completes the proof of the theorem.
  \qed
\end{proof}

For $\Delta=2,4,5,6,7,8,9,\dots$ the following sequence is obtained for the
 coefficients of $n$, which is strictly increasing for every $k \geq 4$:
 $$
  \frac{1}{2} \ < \ \frac{k+1}{2k} \ < \ \frac{k+3}{2k+2} \ < \
     \frac{2k}{3k-1} \ < \ \frac{2k+2}{3k+1} \ < \
       \frac{3k-1}{4k-2} \ < \ \frac{3k+1}{4k} \ < \ \cdots
 $$

\medskip

%
%

To derive a general lower bound in terms of minimum and maximum degree,
we need the following extremal result due to Erd\H os and Gallai.

\begin{theorem}[\cite{erga59}]
\label{erga-thm}
If\/ $G$ is a graph on\/ $n$ vertices that does not contain
a path of order\/ $k$, then it cannot have more than\/ $\frac{n(k-2)}{2}$ edges. Moreover,
the bound is achieved when the graph consists of vertex-disjoint cliques on\/ $k-1$ vertices.
\end{theorem}

Now we are in a position to extend the inequality (\ref{eq:reglow}) from
 regular graphs to arbitrary graphs, using the approach in \cite{Br-2013}.

\begin{theorem}
Let\/ $G$ be a graph of order\/ $n$ with minimum degree\/ $\delta$ and maximum degree\/ $\Delta$.
If\/ $k\ge 3$ and\/ $\delta \ge k-1$, then
$$\psi_k(G) \ge \frac{\delta-k+2}{\delta+\Delta-k+2}\; n.$$
\end{theorem}

\begin{proof}
Let $T \subseteq V(G)$ be a minimum $k$-path vertex cover of $G$, i.e.\ $|T|=\psi_k(G)$, and let $\overline{T}=V(G)-T$.
Further, let $E_T$ and $E_{\overline{T}}$ be the set of edges
 induced by $T$ and by $\overline{T}$,
respectively, and $E_{T\overline{T}}$ be the set of edges with one endvertex in $T$ and the other in $\overline{T}$.
We observe that
$$\sum_{v \in V(G)}d(v)=2|E(G)|=2|E_T|+2|E_{T\overline{T}}|+2|E_{\overline{T}}|$$
and also that $\Delta\cdot |T| \geq 2|E_T|+|E_{T\overline{T}}|$
and $\delta\cdot |\overline{T}| \leq |E_{T\overline{T}}|+2|E_{\overline{T}}|$.
Since the subgraph induced by $\overline{T}$ with edge set $E_{\overline{T}}$ does not contain a path
$P_k$,  Theorem~\ref{erga-thm} implies $|E_{\overline{T}}| \leq \frac{|\overline{T}|(k-2)}{2}$.
From these formulas we get
\begin{align*}
|T| &\geq \frac{1}{\Delta}\left(2|E_T|+|E_{T\overline{T}}|\right) \geq \frac{1}{\Delta}|E_{T\overline{T}}|
    \geq \frac{1}{\Delta}\left(\delta |\overline{T}|-2|E_{\overline{T}}|\right)\\
		&\geq \frac{1}{\Delta}\left(\delta |\overline{T}|-|\overline{T}|(k-2)\right) = \frac{\delta -k+2}{\Delta}\,|\overline{T}|.
\end{align*}
It follows that
$$n =
|T|+|\overline{T}|\le |T|+\frac{\Delta}{\delta -k+2}\, |T|=\frac{\delta + \Delta -k+2}{\delta -k+2}\, |T|$$
and
$$\psi_k(G) = |T| \geq \frac{\delta -k+2}{\delta + \Delta -k+2}\; n.$$
\hfill $\Box $
\bigskip
\end{proof}

\subsection{In terms of degree sequence and average degree}

In this subsection, we improve the following upper bound that was proved in \cite{Br-2011}:
	$$\psi_k(G) \le n- \frac{2k-2}{k}  \sum_{v\in V(G)}\frac{1}{1+d(v)}.$$

\begin{theorem}
	Let\/ $G$ be a graph of order\/ $n$, without isolated vertices.
	 Then for every integer\/ $k\ge 3$ we have
	$$\psi_k(G) \le n- \frac{2k-3}{k-1} \sum_{v\in V(G)}\frac{1}{1+d(v)}.$$
\end{theorem}

\begin{proof}
	We take an ordering $v_1,v_2,\dots,v_n$ of the vertex set at random,
	each of the $n!$ orders being equally likely.
	Having an ordering, we define a weight function\break
	 $w: V(G) \rightarrow \{0, \frac{k-2}{k-1}, 1\}$ as follows. For each vertex $v_i$, let $w(v_i)=1$  if   $v_i$   has no neighbor among $v_1, \dots, v_{i-1}$; let $w(v_i)=(k-2)/(k-1)$ if $v_i$ has exactly one neighbor with an index smaller than $i$; and let $w(v_i)=0$ otherwise.
	Define $Y$ as the set of vertices assigned with positive weights.
	Hence, a vertex $v$ belongs to $Y$ if and only if it appears as first or second vertex 	when the order is restricted to the closed neighborhood $N[v]$ of $v$.
	Since every order has the same probability, also the permutation of $N[v]$
	is a random one and thus, each of the $d(v)+1$ positions is
	equally likely to be taken by $v$.
	Consequently the expected total weight  is
	$$E=\left(1+\frac{k-2}{k-1}\right)\sum_{v\in V(G)} \frac{1}{1+d(v)}.$$
	Expectation of a random variable means a certain average over all outcomes of the
	event, thus there must exist a vertex order under which the total weight on the vertices of $Y$ is at least $E$.
	Observe that $G[Y]$ is a forest.
	Indeed, should a cycle occur, its vertex of largest index would have at least two
	neighbors preceding it, a contradiction to the definition of $Y$. Remark further that any component of $G[Y]$ may contain at most one vertex of weight $1$.  Otherwise adding an edge between two such vertices would result in the same set $Y$ and the subgraph induced by $Y$ would contain a cycle which is a contradiction.
	
	Now, we prove that $G[Y]$ contains a set $S$ of vertices such that $G[S]$ is $P_k$-free and $|S| \ge \sum_{v\in Y} w(v)$. 
	It suffices to prove the analogous statement for the components of $G[Y]$. 
	If $G[Y_i]$ is a component that does not contain any $P_k$-subgraphs,
	 we simply set $S_i=Y_i$ and $|S_i| \ge  \sum_{v\in Y_i} w(v)$ clearly holds. 
	Otherwise, consider a minimum $k$-path vertex cover $T_i$ and its complement $S_i$ in the tree $G[Y_i]$. By Theorem~\ref{p:tree}, we have $|S_i|\ge (k-1)|T_i|$ and under the present assumption $|S_i| \ge k-1$ holds.  Moreover, at most one vertex in $Y_i$ may have a weight of $1$. We infer the following:
		\begin{align*}
	\sum_{v\in Y_i}w(v) & = \sum_{v\in T_i}w(v) + \sum_{v\in S_i}w(v) \le
	  \frac{k-2}{k-1} 
	    \, |T_i|
	   + \frac{k-2}{k-1} 
	    \, |S_i| + \frac{1}{k-1} \\
	& \le \left(\frac{k-2}{(k-1)^2}+ \frac{k-2}{k-1}\right)\, |S_i| + \frac{1}{k-1} = |S_i|- \frac{1}{(k-1)^2}\,|S_i|  + \frac{1}{k-1} \le |S_i|.	  
	\end{align*}
	The set $S=\bigcup S_i$ induces a  $P_k$-free subgraph in $G$ and satisfies 
	$$|S| \ge \sum_{v \in V(G)} w(v) \ge E \ge \frac{2k-3}{k-1}  \sum_{v\in V(G)}\frac{1}{1+d(v)}.$$
	 Since the complement $T= V(G) \setminus S$ is a $k$-path vertex cover in $G$, we conclude that
	 $$\psi_k(G) \le |T| \le  n- \frac{2k-3}{k-1}  \sum_{v\in V(G)}\frac{1}{1+d(v)}.$$
	\qed
	\end{proof}

\begin{corollary}
 If\/ $G$  is an isolate-free graph of order\/ $n$ and average degree\/ $\overline{d}$, then
  $$\psi_k(G)\leq
   \left(1 - \frac{1}{\overline{d}+1}\cdot\frac{2k-3}{k-1}\right) n$$
  for every\/ $k\geq 3$.
\end{corollary}

\section{Chordal graphs}
\label{sec:4}

By definition, a graph is chordal if it does not contain any induced cycles on more than three vertices. In this section we consider the well-studied class of chordal graphs and
 prove upper bounds on $\psi_k$ for its members, in terms of the chromatic number and the order.

\begin{theorem}   \label{t:chi}
If\/ $k\ge 2$ and\/ $G$ is a chordal graph of order\/ $n$ and chromatic number\/ $\chi$, then
$$\psi_k(G) \le \left(1- \frac{2}{\chi}\cdot \frac{k-1}{k}\right) n.$$
Moreover, this bound is tight for every\/ $k$ with\/ $\chi=2$ and also for every\/ $\chi$ with\/ $k=2$.
\end{theorem}

\begin{proof}
Consider a $\chi$-coloring of $G$, and let $X_1$ and $X_2$ be the two largest color classes. Clearly, $|X_1|+|X_2|\ge \frac{2}{\chi}\, n$.
The subgraph induced by $X_1\cup X_2$ is bipartite and chordal. So, $G[X_1\cup X_2]$ is a forest, and it has a $k$-path vertex cover $T$ with $|T|\le (|X_1|+|X_2|)/k$.
Since $(V(G)\setminus (X_1\cup X_2))\cup T$ is a $k$-path vertex cover in $G$, we have
\begin{align*}
\psi_k(G) &\le n-(|X_1|+|X_2|)+\frac{1}{k}\, (|X_1|+|X_2|)\\
				  &\le n-\frac{k-1}{k}\, (|X_1|+|X_2|)\\
					&\le \left(1- \frac{2}{\chi}\cdot \frac{k-1}{k}\right) n.
\end{align*}
Let now $n$ be a multiple of $k\,\chi$.
If $G=P_n$ (hence $\chi=2$), then the equality $\psi_k(G)=n/k$ holds for any $k$ because a vertex has to be chosen into a $k$-path
 vertex cover from any $k$ consecutive vertices of $P_n$.
Also, for $k=2$ with arbitrary $\chi$ let $G=P_n^{\chi-1}$, the graph
 obtained from $P_n$ by joining any two vertices whose
 distance along the path is less than $\chi$.
Then $\psi_2(G)=n-n/\chi$ because $\chi-1$ vertices have to be chosen into
 a $2$-path vertex cover from any $\chi$ consecutive vertices of $P_n^{\chi-1}$.
\hfill $\Box $
\end{proof}

\bigskip

The above inequality is not tight in general; an improvement for almost all
 pairs $k,\chi$ (but not for all, cf.\ Table \ref{tab:comp})
  is given by the following result.
We state it in terms of the clique number, due to the nature of its proof.
But every chordal graph $G$ satisfies $\chi(G)=\omega(G)$, so the two
 theorems might be formulated in the same way.

\begin{theorem}   \label{t:omega}
If\/ $k\ge 3$ and\/ $G$ is a chordal graph on\/ $n$ vertices
 with clique number\/ $\omega$, then
$$\psi_k(G) \le \frac{\omega}{\omega+k-1}\; n.$$
\end{theorem}

\begin{proof}
The validity of the assertion is easy to see if $n\le\omega+k-1$ because
 deleting $n-k+1$ vertices no room remains for $P_k$,
 hence $\psi_k(G) \le n-k+1 \le \frac{\omega}{\omega+k-1}\, n$ if $n$ is so small.
In the rest of the proof we assume $n>\omega+k-1$ and apply induction on $n$.

The following argument is inspired by the theory of treewidth, and in particular
 by the notion of nice tree decomposition, but the formalism here is simpler
 because the graph $G$ is assumed to be chordal.
A detailed general treatment of tree decompositions can be found e.g.\ in the
 monograph \cite{Kloks}, therefore we shall adopt facts from it without
 further references.

Let $v_1,\dots,v_n$ be the vertices of $G$.
Then there exists a rooted binary tree, say $F$, and subtrees $F_1,\dots,F_n$
 in it, where subtree $F_i$ will represent vertex $v_i$ for $i=1,\dots,n$,
 with the following properties.
For every node $x$ of $F$, let us denote $S_x=\{i\mid x\in F_i\}$.
 \begin{itemize}
  \item $V(F_i)\cap V(F_j)\neq\emptyset$ if and only if $v_iv_j\in E(G)$.
  \item $|S_x|\le\omega$ for all $x\in V(F)$.
  \item If $x$ has two children $x'$ and $x''$, then $S_x=S_{x'}=S_{x''}$.
  \item If $x$ has one child $x'$, then $S_x\subset S_{x'}$ or
   $S_{x'}\subset S_x$, moreover $|S_x|=|S_{x'}|+1$ or $|S_x|=|S_{x'}|-1$.
 \end{itemize}

For every $x\in V(F)$ let $V_x$ denote the set of those vertices $v_i$
 for which there exists an $x'$ in the subtree rooted at $x$ in $F$ such that
 $i\in S_{x'}$.
Then we have, in particular, that $|V_x|=n>\omega+k-1$ if $x$ is the root of $F$,
 while $|V_x|=|S_x|\le\omega$ if $x$ is a leaf.
By the property expressed in the fourth bullet above, $|V_x|$ can increase
 by at most 1 when we move from a node to its parent towards the root.
Thus we can choose a node $x$ such that $|V_x|=\omega+k-1$.
Now we put all $v\in S_x$ into the $k$-path vertex cover $T$ to be constructed,
 delete the entire $V_x$ from $G$, moreover omit those $F_i$ which belong to
 vertices $v_i\in V_x$.
After this, we keep only those nodes of $F$ which are incident with at least one
 of the remaining subtrees.

Observe that the $F_j$ not containing $x$ either are entirely in the subtree
 rooted at $x$ (hence they are omitted during this step) or are disjoint
 from the subtree rooted at $x$.
Consequently the vertex set corresponding to $S_x$ is a vertex cut in $G$,
 and there is no edge from $V_x\setminus T$ to $G-V_x$.
It follows that $G[V_x\setminus T]$ is $P_k$-free, and every $k$-path vertex cover
 of $G-V_x$ is completed to a $k$-path vertex cover of $G$ with the current
 vertices of $T$.
Thus, applying the induction hypothesis we obtain
 $$
   \psi_k(G) \le |S_x| + \frac{\omega}{\omega+k-1}\, (n-|V_x|) \le
    \frac{\omega}{\omega+k-1}\, n.
 $$
\hfill $\Box$
\end{proof}

\begin{center}
\noindent
\begin{table}
\begin{center}
\begin{tabular}{|c||cccc|cccc|}
\hline 
$\chi=\omega$ & 2 & 3 & 4 & 5 & 2 & 3 & 4 & 5 \\ 
\hline 
\hline 
$k=2$ & \textbf{1/2} & \textbf{2/3} & \textbf{3/4} & \textbf{4/5} & 2/3 & 3/4 & 4/5 & 5/6 \\ 
\hline 
$k=3$ & \textbf{1/3} & \textbf{5/9} & \textbf{2/3} & 11/15 & 1/2 & 3/5 & \textbf{2/3} & \textbf{5/7} \\ 
\hline 
$k=4$ & \textbf{1/4} & \textbf{1/2} & 5/8 & 7/10 & 2/5 & \textbf{1/2} & \textbf{4/7} & \textbf{5/8} \\ 
\hline 
$k=5$ & \textbf{1/5} & 7/15 & 3/5 & 17/25 & 1/3 & \textbf{3/7} & \textbf{1/2} & \textbf{5/9} \\ 
\hline 
$k=6$ & \textbf{1/6} & 4/9 & 7/12 & 2/3 & 2/7 & \textbf{3/8} & \textbf{4/9} & \textbf{1/2} \\ 
\hline 
\end{tabular} 
\end{center}
\caption{Comparison of Theorem \ref{t:chi} and Theorem \ref{t:omega},
 indicating best current bounds \label{tab:comp}}
\end{table}
\end{center}

\section{Concluding remarks and open problems}
   \label{sec:5}

In this concluding section we discuss a method which leads to further improved upper bounds, study the case of planar graphs, and raise some problems which remain open for future research.

\subsection{Bounds from maximum degree}

An interesting aspect of Theorem \ref{delta-thm} is that it offers a self-improving scheme.
We illustrate this with the case of $\Delta=11$.

\begin{theorem}
If\/ $G$ is a graph on\/ $n$ vertices and\/ $\Delta(G)=11$, then
  $$
    \psi_k(G) \leq \frac{3k+5}{4k+4} \,n
  $$
 for every\/ $k \geq 6$.
\end{theorem}

\begin{proof}
By Theorem \ref{Lovasz-decomp}, $G$ admits a vertex partition $A\cup B$ such that
 both $G[A]$ and $G[B]$ have maximum degree at most 5.
Assume without loss of generality that $|A|\leq n/2\leq |B|$ holds, and let $T_B$ be a
 $k$-path vertex cover in $G[B]$.
By Theorem \ref{delta-thm} we can ensure
 $$
   |T_B|\leq \frac{k+3}{2k+2}\,|B| = |B| - \frac{k-1}{2k+2}\,|B|\leq |B| - \frac{k-1}{4k+4}\,n.
 $$
Since $A\cup T_B$ is a $k$-path vertex cover, and $|A|+|B|=n$,
 the claimed upper bound follows.
\qed
\end{proof}

Certainly the inequality in the theorem is valid also for $k\leq 5$, but to beat the earlier
 formula $\frac{4k}{5k-1}\,n$ we need $(k-1)(k-5)>0$.

It is not clear at the moment, how strong upper bound can be derived as a limit of this
 approach as $\Delta$ gets large, and how restricted it will be as regards the value of $k$.

We should also note that it was necessary to exclude $\Delta=3$ from the range of validity
 in Theorem \ref{delta-thm}.
Indeed, the formula in $(ii)$ would then yield $\frac{4}{k+3}\,n$, i.e., the coefficient
 of $n$ would tend to zero as
 $k$ gets large.
But this would not be a valid formula, as the next result shows.

\begin{theorem}
 For every\/ $k$, there are infinitely many\/ $3$-regular connected graphs\/ $G$ such that\/
  $\psi_k(G) > \frac{1}{4}\,|V(G)|$.
\end{theorem}

\begin{proof}
 Fix any $k\geq 3$.
It follows from results of Erd\H os and Sachs \cite{ersa63} that for every
 even $n\geq 2^{k}$ there exists a 3-regular graph $G$ of order $n$
 such that the length of the shortest cycle (i.e., the girth of $G$) is larger than $k$.
Let now $T$ be any $k$-path vertex cover in $G$, and consider the induced subgraph $H=G-T$.
Say, $H$ has connected components $H_1,\dots,H_s$.
Each $H_i$ is a tree (of diameter at most $k-2$).
Regularity of degree 3 implies that there are exactly $|V(H_i)|+2$ edges
 --- that is, more than the order of $H_i$ --- from
 $V(H_i)$ to $V(G)\setminus V(H_i)$.
Altogether more than $|V(H)|$ edges go from $H$ to the vertices of $T$, while at most $3|T|$ edges go from $T$ to $H$.
It follows that $|T|>\frac{1}{3}\,|V(H)|$ and hence $|T|>\frac{1}{4}\,n$.
\qed
\end{proof}

\subsection{Planar graphs} \label{sec:planar}

In this subsection we present some partial results related to a conjecture which was first mentioned by Bre\v sar et al.\ in 2011. Note that it is easy to construct planar graphs that satisfy the relation with equality, e.g.\ the vertex-disjoint union of octahedron graphs, adjacent in any planar way.
\begin{conjecture}[\cite{Br-2011}]
	\label{conj:planar}
	If\/ $G$ is a planar graph of order\/ $n$, then\/
	$\psi_3(G) \le \frac{2}{3}\, n$.
\end{conjecture}

It is a basic fact that $\omega(G) \le 4$ holds for every planar graph. Then, Theorem~\ref{t:omega} allows us to conclude that Conjecture~\ref{conj:planar} is true for those planar graphs which are also chordal. We may also derive upper bounds on $\psi_k(G)$ for $k \ge 4$.

\begin{corollary}
	If\/ $G$ is a graph of order\/ $n$ which is chordal and planar, then
	$$\psi_3 (G) \le \frac{2}{3}\; n , \qquad
	\psi_4 (G) \le \frac{4}{7}\; n ,
	\quad {\rm and} \quad \psi_5 (G) \le \frac{1}{2}\; n .$$
\end{corollary}

\medskip
On the class of $K_3$-free planar graphs Theorem~\ref{n-m-thm} and the inequality $m \leq 2n-4$ together imply 
$$\psi_3(G) \leq \frac{4n+m}{9} \leq \frac{4n+2n-4}{9}=\frac{2}{3}n-\frac{4}{9}.$$
Therefore, Conjecture~\ref{conj:planar} is also true for $K_3$-free planar graphs. Moreover, Corollary~\ref{cor:31}  will show that the upper bound $2n/3$ can be considerably improved on this graph class.

\medskip

The \emph{forest number} $a(G)$ of a graph G is defined as the maximum number of vertices in an induced forest of $G$. As it is mentioned shortly in the concluding section of \cite{Br-2011}, there is an important relation between the $k$-path vertex cover number and the forest number. For the sake of completeness, we present a proof for this statement.  

\begin{proposition}[\cite{Br-2011}]  \label{forest-number}
	For every integer\/ $k\ge 3$ and for every graph\/ $G$ of order\/ $n$,
	$$\psi_k(G) \le n- \frac{k-1}{k}\; a(G).$$
\end{proposition}

\begin{proof} Consider a set $S$ of vertices such that $|S|=a(G)$ and $G[S]$ is a forest. By Proposition~\ref{p:tree}, every $k$-path inside $G[S]$ can be covered by a set $T$ of at most $|S|/k$ vertices. The set $(V(G)\setminus S) \cup T$ is clearly a $k$-path vertex cover in $G$. Hence, we have
	$$\psi_k(G) \le n-|S| + \frac{|S|}{k} = n- \frac{k-1}{k} |S|=
	n - \frac{k-1}{k}\; a(G).$$
	\qed
\end{proof}

The following is a famous unsolved problem on the forest number of planar graphs.
\medskip

\noindent {\bf Albertson--Berman Conjecture (\cite{AB})} {\it If\/ $G$ is
	a planar graph of order\/ $n$, then\/ $a(G) \ge n/2$.}

\bigskip
Proposition~\ref{forest-number} shows that the Albertson--Berman Conjecture, if true, implies Conjecture~\ref{conj:planar}. Combining Proposition~\ref{forest-number} and the following two results we may also conclude some upper bounds on $\psi_k(G)$.
\begin{theorem}[\cite{Borodin}]
	\label{borodin}
	If\/ $G$ is a planar graph of order\/ $n$, then\/ $a(G) \ge \frac{2}{5}\,n$.
\end{theorem}

\begin{theorem}[\cite{KLS}]
	\label{KLS}
	If\/ $G$ is a planar graph of order\/ $n$ that is\/ $K_3$-free, then\/ $a(G) \ge
	\frac{71}{128}\, n + \frac{9}{16}$.
\end{theorem}


\begin{corollary} If\/ $G$ is a planar graph of order\/ $n$, then
	$$\psi_3 (G) \le \frac{11}{15}\; n< 0.734\, n$$
	and
	$$\psi_6 (G) \le \frac{2}{3}\; n.$$
\end{corollary}

\begin{corollary} \label{cor:31}
	If\/ $G$ is a planar and\/ $K_3$-free graph on\/ $n$ vertices, then
	$$\psi_3 (G) \le \frac{121}{192}\; n-\frac{3}{8}.$$
\end{corollary}

\bigskip

\subsection{Further conjectures and open problems} \label{sec:chordal-planar}

Some of the estimates proved in this paper are tight, while some others probably aren't.
In particular, on the class of chordal graphs we expect that the following improvement can be made.

\begin{conjecture}
If\/ $k\ge 3$ and\/ $G$ is a chordal graph with clique number\/ $\omega$ and $n=|V(G)|$, then
$$\psi_k(G) \le \frac{\omega-1}{\omega+k-2}\; n.$$
\end{conjecture}

The analogous inequality for $k=2$ can easily be seen for every graph in which the equality $\chi=\omega$ is valid --- hence for all perfect graphs and in particular for every chordal graph --- by taking a proper vertex coloring with $\omega$ colors and deleting the largest color class.

The following further problems arise in a natural way in connection with our results on planar graphs.

\begin{problem}
 Determine the smallest constants\/ $c_{1,k}$ and\/ $c_{2,k}$ such that
  \begin{itemize}
   \item[$(i)$] if\/ $G$ is planar, then\/ $\psi_k (G) \le (c_{1,k}+o(1))\, n$;
   \item[$(ii)$] if\/ $G$ is planar and\/ $K_3$-free, then\/
    $\psi_k (G) \le (c_{2,k}+o(1))\, n$,
  \end{itemize}
 as $n\to\infty$.
\end{problem}

Note that Conjecture~\ref{conj:planar} formulates the explicit guess $c_{1,3}=2/3$.
\medskip

Corollary~\ref{many-bounds} offers various upper bounds on $\psi_3(G)$. The following problem naturally arises.
\begin{problem}
Are the upper bounds stated in Corollary~\ref{many-bounds} sharp? 
\end{problem}


\end{document}